\def\timestamp{%
Time-stamp: <span.tex: Thursday 11-03-2004 at 11:57:18 (cet)>}
\def\stripname Time-stamp: <#1 #2>{#2}
\edef\thefiledate{\expandafter\stripname\timestamp}
\documentclass[a4paper]{amsart}
\usepackage{amsfonts,amscd}

\theoremstyle{plain}
\newtheorem{theorem}{Theorem}[section]
\newtheorem{lemma}[theorem]{Lemma}

\newtheorem{corollary}[theorem]{Corollary}

\theoremstyle{definition}
\newtheorem{question}[theorem]{Question}

\newcounter{claim}[theorem]
\renewcommand\theclaim{\arabic{claim}}
\newenvironment{claim}%
           {\endgraf\refstepcounter{claim}\smallskip
            \def\proof{\endgraf\vskip1pt plus1pt \noindent\textsl{Proof}. }
            \noindent\textsl{Claim \theclaim}. \ignorespaces}%
           {\parfillskip0pt\hfil$\symmdif$\endgraf\smallskip}

\newcommand{\powerinfront}[2]{{}^#1 #2}
\newcommand{\interior}{\operatorname{int}}
\newcommand{\cl}{\operatorname{cl}}

\let\le\undefined
\let\ge\undefined

\DeclareMathSymbol\restr{2}{AMSa}{"16}
\DeclareMathSymbol\symmdif{2}{AMSa}{"4D}
\DeclareMathSymbol\le{3}{AMSa}{"36}
\DeclareMathSymbol\ge{3}{AMSa}{"3E}

\DeclareMathSymbol\CC{0}{AMSb}{`C}
\DeclareMathSymbol\HH{0}{AMSb}{`H}
\DeclareMathSymbol\II{0}{AMSb}{`I}
\DeclareMathSymbol\MM{0}{AMSb}{`M}
\DeclareMathSymbol\reals{0}{AMSb}{`R}
\newcommand\Hstar{\HH^*}

\newcommand\cB{\mathcal{B}}
\newcommand\cC{\mathcal{C}}
\newcommand\cR{\mathcal{R}}

\newcommand\CH{\mathsf{CH}}
\newcommand\ZFC{\mathsf{ZFC}}

\usepackage[numeric]{amsrefs}[2002/02/28]

\begin{document}

\title{Span, chainability and the continua $\Hstar$ and~$\II_u$}
\author{K. P. Hart}
\author{B. J. van der Steeg}
\address{Faculty EEMCS\\
         TU Delft\\
         Postbus 5031\\
         2600~GA {} Delft\\
         the Netherlands}

\email[K. P. Hart]{k.p.hart@ewi.tudelft.nl}
\email[B. J. van der Steeg]{berd@punch-basketball.nl}

\urladdr[K. P. Hart]{http://aw.twi.tudelft.nl/\~{}hart}

\dedicatory{Dedicated to the memory of Zoli Balogh}

\begin{abstract}
We show that the continua $\II_u$ and $\Hstar$ are
non-chainable and have span nonzero.
Under $\CH$ this can be strengthened to surjective symmetric
span nonzero.

We discuss the logical consequences of this.
\end{abstract}



\date{\thefiledate}

\maketitle

\section{Introduction}\label{sect-intro}

Chainable (or arc-like) continua are `long and thin'; in an attempt to capture
this idea in metric terms Lelek introduced, in~\cite{L1}, the notion of span.
Chainable continua have span zero, which is useful in proving that certain
continua are not chainable.
The converse, a conjecture by Lelek in~\cite{L2}, is one of the 
main open problems in continuum theory today.
While the particular value of the span of a continuum depends on the metric
chosen, the distinction between span zero and span nonzero is a topological 
one.
As chainability is a topological notion as well, Lelek's theorem and conjecture
are meaningful in the class of all Hausdorff continua.
We investigate the chainability and span of several continua that are closely
connected to the \v{C}ech-Stone compactification of the real line.

\section{Preliminaries}\label{sect-prelim}

\subsection{Various kinds of span}

The kinds of span that we consider in this paper are, in the metric case,
defined as suprema of distances between the diagonal of the continuum
and certain subcontinua of the square.
The following families of subcontinua feature in these definitions:
\begin{description}
  \item[$S(X)$] the symmetric subcontinua of~$X^2$, i.e., those that
                satisfy $Z=Z^{-1}$;
  \item[$\Sigma(X)$] the subcontinua of~$X^2$ 
                     that satisfy $\pi_1[Z]=\pi_2[Z]$; and
  \item[$\Sigma_0(X)$] the subcontinua of~$X^2$ 
                     that satisfy $\pi_2[Z]\subseteq\pi_1[Z]$.
\end{description}
Here, $\pi_1$ and $\pi_2$ are the projections onto the first and second
coordinates respectively.
It is clear that $S(X)\subseteq\Sigma(X)\subseteq\Sigma_0(X)$ and hence
that $s(X)\le\sigma(X)\le\sigma_0(X)$, where
\begin{enumerate}
  \item $s(X)=\sup\bigl\{d(\Delta(X),Z):Z\in S(X)\bigr\}$;
  \item $\sigma(X)=\sup\bigl\{d(\Delta(X),Z):Z\in \Sigma(X)\bigr\}$; and
  \item $\sigma_0(X)=\sup\bigl\{d(\Delta(X),Z):Z\in \Sigma_0(X)\bigr\}$.
\end{enumerate}
These numbers are, respectively, the \emph{symmetric span},
the \emph{span} and the \emph{semi-span} of~$X$.

If one uses, in each definition, only the continua~$Z$ with $\pi_1[Z]=X$
then one gets the \emph{surjective symmetric span}, $s^*(X)$,
the \emph{surjective span}, $\sigma^*(X)$, 
and the \emph{surjective semi-span}, $\sigma_0^*(X)$, of~$X$
respectively.
The following diagram shows the obvious relationships between the six
kinds of span.
\begin{equation}\label{eq:1}
\begin{CD}
  s(X)   @>>> \sigma(X)   @>>> \sigma_0(X)   \\
  @AAA        @AAA             @AAA          \\
  s^*(X) @>>> \sigma^*(X) @>>> \sigma_0^*(X) \\
\end{CD}
\end{equation}
Topologically we can only distinguish between a span being zero or nonzero.
A span is zero if and only if every continuum from its defining family
intersects the diagonal.
This defines span zero (or span nonzero) for the six possible types of 
span in general continua.

Below we will show that for the continua $\Hstar$ and $\II_u$ all six
kinds of span are nonzero.
Diagram~(\ref{eq:1}) shows that it will be most difficult to show that
$s^*$ is nonzero (or dually that it would be hardest to show that
$\sigma_0$ is zero).
Indeed, we will give successively more difficult proofs that the various spans
are nonzero, where we traverse the diagram from top right to bottom left.

The need for these different proofs lies in their set-theoretic assumptions.
We need nothing beyond $\ZFC$ to show that $\sigma^*(\Hstar)$ and
$\sigma(\II_u)$ are nonzero; to show that the other spans (in particular $s^*$)
are nonzero we shall need the Continuum Hypothesis ($\CH$). 

\subsection{Chainability}

A continuum is \emph{chainable} if every open cover of it has an open
refinement that is a chain cover, where 
$\cC=\{C_1,\ldots,C_m\}$ \emph{chain cover} 
if $C_i\cap C_j$ is nonempty if and only if $|i-j|\le 1$.

One readily shows that every chainable continuum has span zero, whatever
kind of span one uses.
This follows from the fact that chainability is a hereditary property
of continua and from the following theorem whose proof we give for
completeness sake.

\begin{theorem}\label{thm:chainable.spanzero}
Every chainable continuum has surjective semi-span zero.
\end{theorem}
\begin{proof}
Let $X$ be a chainable continuum and let $Z$ be a subcontinuum of~$X^2$
that is disjoint from $\Delta(X)$.
Let $\mathcal{U}$ be a finite open cover of $X$ such that 
$U^2\cap Z=\emptyset$ for all $U\in\mathcal{U}$.
Next let $\{V_1,V_2,\ldots,V_n\}$ be an open chain cover that 
refines~$\mathcal{U}$.
Define open sets $O_1$ and $O_2$ in $X^2$ by
$$
O_1=\bigcup\{V_i\times V_j:i<j\},\quad O_2=\bigcup\{V_i\times
V_j:i>j\}.
$$
Then $Z\subset O_1\cup O_2$ and $Z\cap O_1\cap O_2=\emptyset$.
As $Z$ is connected, it is contained in one of $O_1$ or $O_2$,
say $Z\subseteq O_2$. 
Then $\pi_1[Z]\subseteq \bigcup_{i<n}V_i$ and
     $\pi_2[Z]\subseteq \bigcup_{i>1}V_i$.
This means that neither $\pi_1[Z]$ nor $\pi_2[Z]$ is equal to~$X$.
\end{proof}

\subsection{The continua $\II_u$ and $\Hstar$}

In this paper we will be investigating the different kinds of span
and the chainability of the continua $\II_u$ and~$\Hstar$. 
These two spaces are related to one another.
Following~\cite{Mio} and~\cite{Ha}, we will use the space
$\MM=\omega\times\II$ in our investigation of the
spaces $\II_u$ and~$\Hstar$, where $\II$
denotes the unit interval $[0,1]$.

The map $\pi:\MM\rightarrow\omega$ given by $\pi(n,x)=n$ is
perfect and monotone, as is its \v{C}ech-Stone extension $\beta\pi$. 
The preimage of an ultrafilter $u\in\beta\omega^*$ is a continuum
and denoted by $\II_u$.

Given any sequence $\langle x_n\rangle_{n\in \omega}$ in $\II$ and any
$u\in\omega^*$ there is a unique point, denoted~$x_u$, in~$\II_u$ such
that for every $\beta\MM$-neighborhood~$O$ of~$x_u$, the set
$\{n\in\omega:(n,x_n)\in O\}$ is an element of~$u$, i.e.,
$x_u$~is the $u$-limit of the sequence $\langle (n,x_n)\rangle_{n\in \omega}$. 
These points form a dense set $\CC_u$ of cut~points of~$\II_u$, 
for details see~\cite{Ha}.
The set~$\CC_u$ is in fact the ultrapower of~$\II$ by the 
ultrafilter~$u$, i.e., the set $\powerinfront\omega\II$ modulo the equivalence
relation $x\sim_uy$ defined by $\{n:x_n=y_n\}\in u$.

The continuum $\II_u$ is irreducible between the points~$0_u$ and~$1_u$
(defined in the obvious way) and as it has a natural pre-order $\le_u$
defined by $x\le_u y$ iff every subcontinuum of $\II_u$ that contains
$0_u$ and $y$ also contains~$x$.
The equivalence classes under the equivalence relation 
``$x\le_uy$ and $y\le_ux$'' are called layers and the set of layers is linearly
ordered by~$\le_u$.
The points of $\CC_u$ provide one-point layers, the restriction 
of~$\le_u$ to this set coincides with the ultrapower order defined
by $\{n:x_n\le y_n\}\in u$.
We shall freely use interval notation, allowing non-trivial layers as 
end~points.

If $\langle x_n\rangle_{n\in\omega}$ is a strictly increasing sequence
in~$\II_u$ then its supremum~$L$ is a non-trivial layer.
Because $\beta\MM\setminus\MM$ is an $F$-space the closure of 
$\{x_n:n\in\omega\}$ is homeomorphic to~$\beta\omega$; by upper semicontinuity
the remainder (which is a copy of $\omega^*$) must be contained in~$L$.
We call such a layer a countable-cofinality layer.

The continuum $\Hstar$ is the remainder of the
\v{C}ech-Stone compactification $\beta\HH$, where
$\HH$ is the half line $[0,\infty)$. Let
$q:\MM\rightarrow\HH$ be given by $q(n,x)=n+x$, then
$q$ is a perfect map and its \v{C}ech-Stone extension $\beta
q:\beta\MM\rightarrow\beta\HH$ maps $\MM^*$
onto~$\Hstar$. Again for properties of $\Hstar$ and
its relation to $\II_u$ see~\cite{Ha}.

\section{The span of $\Hstar$}

In this section we show that the surjective (semi-)span of $\Hstar$ is 
nonzero.
The following theorem more than establishes this.

\begin{theorem}\label{thm:hstar-not-fpp}
There exists a fixed-point free autohomeomorphism of\/ $\Hstar$.
\end{theorem}

\begin{proof}
Let $f:\HH\to\HH$ be the map defined by $f:x\mapsto x+1$.
It is clear that $\beta f$ maps $\Hstar$ onto $\Hstar$. 
The restriction $f^*=\beta f\restr\Hstar$ is a fixed-point free
autohomeomorphism of~$\Hstar$.

To see that $f^*$ is an autohomeomorphism consider $g:\HH\to\HH$
defined by $g(x)=\max\{0,x-1\}$.
From the fact that $f\bigl(g(x)\bigr)=x$ and $g\bigl(f(x)\bigr)=x$ 
for $x\ge1$ it follows that
$f^*\circ g^*$ and $g^*\circ f^*$ are the identity on~$\Hstar$.

That $f$ is fixed-point free on~$\Hstar$ follows by considering the following
closed cover $\{F_0,F_1,F_2,F_3\}$ of~$\HH$, defined by
$F_i=\bigcup_n[2n+\frac i2,2n+\frac{i+1}2]$.
Observe that $f^*[F_i^*]=F_{i+2\bmod4}^*$ 
and that $F_i^*\cap F_{i+2\bmod4}^*$ is always empty, so that $f^*(x)\neq x$
for $x\in\Hstar$.
\end{proof}

\begin{corollary}\label{cor:sigmastar.Hstar.nonzero}
$\sigma^*(\Hstar)$ is nonzero.
\end{corollary}

\begin{proof}
The graph of $f^*$ is a continuum in $\Hstar\times\Hstar$ that is disjoint
from the diagonal and whose projection on each of the axes is~$\Hstar$.
\end{proof}

Later we shall see that under $\CH$ even $s^*(\Hstar)$ is nonzero.

By Theorem~\ref{thm:chainable.spanzero} we also know that $\Hstar$
is not chainable.
The reader may enjoy showing that the four open sets $U_0$, $U_1$, $U_2$
and $U_3$ defined by 
$$
U_i=\bigcup_{n<\omega}(8n+2i, 8n+2i+3)
$$
induce an open cover of~$\Hstar$ without a chain refinement.

\subsection{More fixed-point free homeomorphisms}

We use the description of indecomposable subcontinua from~\cite{DH1}
to show that many subcontinua of~$\Hstar$ have fixed-point free
autohomeomorphisms.

We use the shift-map $\sigma:\omega\to\omega$, defined by $\sigma(n)=n+1$,
and its extension to~$\beta\omega$.
We note that $\sigma$~is an autohomeomorphism of~$\omega^*$.
We also write $u+1$ for $\sigma(u)$ and $u-1$ for $\sigma^{-1}(u)$.

For $F\subseteq\omega^*$ we put $\MM_F=\bigcup_{u\in F}\II_u$
and $C_F=\beta q[\MM_F]$.
We say that $F$~is $\sigma$-invariant if $u+1,u-1\in F$ whenever $u\in F$.
Clearly then, if $F$~is $\sigma$-invariant then $f^*\restr C_F$~is an
autohomeomorphism of~$C_F$,
where $f^*$~is the autohomeomorphism of~$\Hstar$ defined in the proof
of Theorem~\ref{thm:hstar-not-fpp}.

From~\cite{DH1} we quote the following:
$C_F$~is a subcontinuum whenever $F$~is closed, $\sigma$-invariant
and not the union of two disjoint proper closed $\sigma$-invariant subsets.
In that case $C_F$ is indecomposable if and only if $F$~is dense-in-itself.

From~\cite{DH1} we also quote:
if $K$~is an indecomposable subcontinuum of~$\Hstar$ then there is a strictly
increasing sequence~$\langle a_n\rangle_n$ in~$\HH$ that diverges to~$\infty$
and such that $K=q_a[C_F]$ for some closed dense-it-itself $\sigma$-invariant
subset~$F$ of~$\omega^*$ that is not the union of two disjoint proper 
closed $\sigma$-invariant subsets and
where $q_a:\Hstar\to\Hstar$ is induced by the piecewise linear self-map
of~$\HH$ that sends $n$ to~$a_n$.

We can combine all this into the following theorem.

\begin{theorem}\label{thm:all.indec.not.fpp}
Every indecomposable subcontinuum of\/~$\Hstar$ has a fixed-point free 
autohomeomorphism (and hence surjective span nonzero).  
\end{theorem}

\section{The span of $\II_u$}

In this section we show that $\II_u$ has span nonzero for any ultrafilter~$u$;
the next section will be devoted to the surjective versions of span.

The following theorem, akin to Theorem~\ref{thm:hstar-not-fpp} and with
a similar proof, provides a continuum witnessing that $\II_u$ has nonzero
span.

\begin{theorem}\label{thm:layer-not-fpp}
Every countable-cofinality layer has a fixed-point free autohomeomorphism.
\end{theorem}

This follows from Theorem~\ref{thm:all.indec.not.fpp} but for later use
we give a direct construction, which establishes a bit more, 
namely that the interval
$[0_u,L]$ has a fixed-point free continuous self-map.

\begin{proof}
We prove the theorem for one particular layer but the argument is easily
adapted to the general case.

For $m\in\omega$ put $x_m=1-2^{-m}$; then $\{x_m\}_{m<\omega}$ is a strictly
increasing sequence in~$\II$ that converges to~$1$ and with $x_0=0$. 
Let $x_{m,u}$ denote the point of $\II_u$ that corresponds to the constant
sequence $\{x_m\}_{n\in\omega}$ in~$\II$.
Then $\{x_{m,u}\}_{m\in\omega}$
is a strictly increasing sequence in $\II_u$; 
let $L$ denote the limit of this sequence, a non-trivial layer of~$\II_u$.

We define a map $f:\II_u\rightarrow\II_u$ by defining it on~$\MM$, 
taking its \v{C}ech-Stone extension and restricting that to~$\II_u$.
\begin{enumerate}
\item Let $f\restr\II_0$ be equal to the identity.
\item For all $n\ge 1$ let $f\restr\II_n$ be the piecewise linear map that
      maps $(n,x_m)$ to $(n,x_{m+1})$ for all $m<n$ and the point $(n,1)$ 
      to itself.
\end{enumerate}

\begin{claim}
The \v{C}ech-Stone extension of the map $f$ maps $[0_u,L]$
homeomorphically onto $[x_{1,u},L]$.
\proof
It is not hard to see that $\beta f$ maps the interval
$[x_{m,u},x_{m+1,u}]$ of $\II_u$ homeomorphically onto $[x_{m+1,u},x_{m+2,u}]$
for all $m\in\omega$. 
This implies that $\beta f$ maps $[0_u,L)$ homeomorphically 
onto $[x_{1,u},L)$. 
The fact that $[0_u,L]=\beta[0_u,L)$ now establishes the claim.
\end{claim}

We let $h$ denote the restriction of $\beta f$ to $[0_u,L]$.
The fact that $[0_u,L]=\beta[0_u,L)$ also establishes the 
following claim.

\begin{claim}
The restriction $h\restr L$ maps $L$ homeomorphically onto $L$.
\end{claim}

To see that $h$ has no fixed points we argue as in the proof of 
Theorem~\ref{thm:hstar-not-fpp}.

For every $m$ let $a_m$ be the mid point of the interval $(x_m,x_{m+1})$.
Note that the map $f$ maps $(n,a_m)$ onto the point $(n,a_{m+1})$
whenever $m<n$. 
Define the following closed subsets $F_i$ for $i=0$, $1$, $2$ and~$3$:
$$\tabskip0pt plus\displaywidth
\halign to\displaywidth{$#$\tabskip0pt&${}#{}$&$#$\hfil&#&
                        \hfil$#$&${}#{}$&$#$\tabskip0pt plus\displaywidth\cr
F_0&=&\bigcup_n\bigl(\{n\}\times\bigcup_{m<n}[x_{2m},a_{2m}]\bigr),&\qquad &
F_2&=&\bigcup_n\bigl(\{n\}\times\bigcup_{m<n}[x_{2m+1},a_{2m+1}]\bigr),\cr
\noalign{\vskip 2pt}
F_1&=&\bigcup_n\bigl(\{n\}\times\bigcup_{m<n}[a_{2m},x_{2m+1}]\bigr),&\qquad &
F_3&=&\bigcup_n\bigl(\{n\}\times\bigcup_{m<n}[a_{2m+1},x_{2m+2}]\bigr).\cr}
$$
Note that the closure in $\beta\MM$ of the union of the $F_i$'s 
contains the interval $[0_u,L]$ of $\II_u$. 
Also note that the closed set $F_i$ is mapped onto the closed set
$F_{i+2\bmod4}$, so $f[F_i]\cap F_i=\emptyset$. 
As in the proof of Theorem~\ref{thm:hstar-not-fpp} this implies
that $h$ has no fixed points.
\end{proof}

As before we get the following corollaries.

\begin{corollary}
The surjective span of $L$ is nonzero, hence $\sigma(\II_u)$ is nonzero.\qed
\end{corollary}

\begin{corollary}
The surjective semi-span of $[0_u,L]$ is nonzero.\qed
\end{corollary}

It will be more difficult to prove the same for $\II_u$.

\section{The surjective spans of $\II_u$ and $\Hstar$}
\label{sec:surj-spans}

Using the map from the previous section and the retraction we get
from the next theorem we will show that under $\CH$ there
exists a fixed-point free continuous self map of $\II_u$; 
as the map is not onto this only implies that the surjective semi-span 
of $\II_u$ is nonzero.
However, the special structure of~$\II_u$ will allow us to build, using
the graph of this map, a symmetric subcontinuum of~$\II_u^2$ that will 
witness $s^*(\II_u)\neq0$;
it will then also be possible to show that $s^*(\Hstar)$ is nonzero.

We retain the notation from the previous section but we write $a_m=x_{m,u}$
for ease of notation and we recall that layer $L$ is the supremum, in $\II_u$,
of the set $\{a_m:m\in\omega\}$.
The following theorem is what makes the rest of this section work.

\begin{theorem}[$\CH$]\label{retraction}
$L$ is a retract of $[L,1_u]$.
\end{theorem}

Before we prove the theorem we give the promised consequences.

\begin{theorem}[$\CH$]
The continuum $\II_u$ does not have the fixed-point property.  
\end{theorem}

\begin{proof}
Let $h:[0_u,L]\to[0_u,L]$ be the map constructed in the proof of 
Theorem~\ref{thm:layer-not-fpp} and let $r:[L,1_u]\to L$ be the retraction
from Theorem~\ref{retraction}.
Extending $r$ by the identity on~$[0_u,L]$ yields a retraction~$r^*$
from~$\II_u$ onto~$[0_u,L]$.
The composition $h\circ r^*$ is then a fixed-point free continuous 
self-map of~$\II_u$.  
\end{proof}

\begin{corollary}[$\CH$]
 The surjective semi-span of $\II_u$ is nonzero.
 \end{corollary}

\begin{proof}
The graph of $h\circ r^*$ is a witness.
\end{proof}

We now show how to make $s^*(\II_u)$ nonzero.

\begin{corollary}[$\CH$]\label{thm:sssI-u.neq.0}
The surjective symmetric span of\/~$\II_u$ is nonzero.  
\end{corollary}

\begin{proof}
Let $G$ be the graph of $h\circ r^*$.
We complete $G$ to symmetric continuum by adding the following continua:
$\{1_u\}\times[0_u,L]$, 
$[h(0_u),1_u]\times\{0_u\}$,
$G^{-1}$,
$[0_u,L]\times\{1_u\}$,
and $\{0_u\}\times[h(0_u),1_u]$.
It is straightforward to check that the union~$Z$ is a continuum
(each continuum meets its successor) that is symmetric and projects
onto each axis.
As none of the pieces intersects the diagonal we get a witness to
$s^*(\II_u)$ being nonzero.
\end{proof}

\begin{corollary}[$\CH$]
The surjective symmetric span of\/~$\Hstar$ is nonzero.  
\end{corollary}

\begin{proof}
We begin by taking the graph $F$ of the map~$f$ from 
Theorem~\ref{thm:hstar-not-fpp} and its inverse $F^{-1}$; unfortunately 
the union $F\cup F^{-1}$ is not connected, as $F$ and $F^{-1}$ are disjoint.
To connect them we take one ultrafilter~$u$ on~$\omega$ and observe that
the image $q[\II_u]$ connects the ultrafilters~$u$ and $u+1$.
The image $K=(q\times q)[Z]$, where $Z$~is from the proof of 
Corollary~\ref{thm:sssI-u.neq.0} meets both~$F$ (in $(u,u+1)$)
and $F^{-1}$ (in $(u+1,u)$).
The union $F\cup K\cup F^{-1}$ is a witness to $s^*(\Hstar)\neq 0$.
\end{proof}

\subsection{Proof of Theorem~\ref{retraction}}

We will construct the retraction by algebraic, rather than topological,
means.
Let $\cR$ be the family of finite unions of closed intervals
of~$\II$ with rational endpoints. 
For every $f\in\powerinfront{\omega}{\cR}$ we define the closed
subset~$A_f$ of~$\MM$ by
$$
A_f=\bigcup_{n<\omega}\{n\}\times f(n).
$$
These sets form a lattice base for the closed sets of~$\MM$,
i.e., it is a base for the closed sets and closed under finite unions
and intersections.
It is an elementary exercise to show that disjoint closed sets in~$\MM$
can be separated by disjoint closed sets of the form~$A_f$.
This implies that the closures $\cl{A_f}$ form a lattice base for the
closed sets of~$\beta\MM$.
It follows that 
$\cB=\{\cl{A_f}\cap L:f\in\powerinfront{\omega}{\cR}\}$
is a base for the closed sets of~$L$ and similarly that
$\cC=\{\cl{A_f}\cap [L,1_u]:f\in\powerinfront{\omega}{\cR}\}$
is a base for~$[L,1_u]$.

Theorem~1.2 from \cite{DH} tells us that in order to construct a retraction
from $[L,1_u]$ onto~$L$ it suffices to construct a map $\varphi:\cB\to\cC$ that
satisfies
\begin{enumerate}
\item $\varphi(\emptyset)=\emptyset$, and
      if $F\neq\emptyset$ then $\varphi(F)\neq\emptyset$;
\item if $F\cup G=L$ then $\varphi(F)\cup\varphi(G)=[L,1_u]$;       
\item if $F_1\cap\cdots\cap F_n=\emptyset$
      then $\varphi(F_1)\cap\cdots\cap\varphi(F_n)=\emptyset$; and
\item $\varphi(F)\cap L= F$.
\end{enumerate}
The retraction $r:[L,1_u]\to L$ is then defined by 
$r(x)={}$`the unique point in $\bigcap\{F:x\in\varphi(F)\}$'.
The first three conditions ensure that $r$~is well-defined, continuous
and onto; the last condition ensures that $r\restr L$ is the identity.

There is a decreasing $\omega_1$-sequence 
$\langle b_\alpha\rangle_{\alpha<\omega_1}$ 
of cut points in~$\II_u$ such that $L=\bigcap_{m,\alpha}[a_m,b_\alpha]$:
by \cite{Ha}*{Lemma 10.1} such a sequence must have uncountable
cofinality and by~$\CH$ the only possible (minimal) length then is~$\omega_1$.
For each $\alpha$ choose a sequence $\langle b_{\alpha,n}\rangle_{n\in\omega}$
in~$\II$ such that $b_\alpha=b_{\alpha,u}$.

Again by $\CH$ we list $\powerinfront{\omega}{\cR}$ in an $\omega_1$-sequence
$\langle f_\alpha\rangle_{\alpha<\omega_1}$.
We will assign to each~$f_\alpha$
a $g_\alpha\in\powerinfront{\omega}{\cR}$ in such a way that
$\cl{A_{f_\alpha}}\cap L \mapsto 
  \cl{A_{g_\alpha}}\cap[L,1_u]$ 
defines the desired map~$\varphi$.

The assignment will be constructed in a recursion of length $\omega_1$,
where at stage $\alpha$ we assume the conditions (1)--(4) are satisfied
for the $A_{f_\beta}$ and $A_{g_\beta}$ with $\beta<\alpha$ and choose
$g_\alpha$ in such a way that they remain satisfied for $\beta\le\alpha$.
At every stage we will list $\alpha$ in an $\omega$-sequence; this means that
it suffices to consider the case $\alpha=\omega$ only.

We need a few lemmas that translate intersection properties in~$\cB$ and~$\cC$
to~$\cR$.

\begin{lemma}\label{lemma:AfcapL=0}
$\cl{A_f}\cap L=\emptyset$ if and only if there are $m$ and $\alpha$
such that the set 
$\bigl\{n:f(n)\cap[a_{m,n},b_{\alpha,n}]=\emptyset\bigr\}$ belongs to~$u$.  
\end{lemma}

\begin{proof}
By compactness $\cl{A_f}\cap L=\emptyset$ if and only if there 
are $m$ and $\alpha$ such that $\cl{A_f}\cap [a_m,b_\alpha]=\emptyset$
and the latter is equivalent
to $\{n:f(n)\cap[a_{m,n},b_{\alpha,n}]=\emptyset\}\in u$
again by compactness and the formula
\begin{equation*}
\cl{A_f}\cap [a_m,b_\alpha]=
 \bigcap_{U\in u}
\cl\Bigl(
   \bigcup_{n\in U}\{n\}\times\bigl(f(n)\cap[a_{m,n},b_{\alpha,n}]\bigl)
   \Bigr).
\qedhere
\end{equation*}
\end{proof}

\begin{lemma}\label{lemma:Af=Ag}
$\cl A_f\cap L=\cl A_g\cap L$ if and only if there are $m$ and $\alpha$ such 
that the set 
$\bigl\{n:f(n)\cap[a_{m,n},b_{\alpha,n}]=g(n)\cap[a_{m,n},b_{\alpha,n}]\bigr\}$  
belongs to~$u$.
\end{lemma}

\begin{proof}
The `if' part is clear.
For the `only if' part let $D$ be the set of all mid points of all maximal
intervals in $A_f\setminus A_g$; then $\cl D\subseteq\cl A_f\setminus\cl A_g$
and so $\cl D\cap L=\emptyset$.
Observe that $D=A_h$ for some $h$, so there are $m$ and $\alpha$
as in Lemma~\ref{lemma:AfcapL=0} for~$D$.
By convexity, for each~$n$ the interval $[a_{m,n},b_{\alpha,n}]$ meets 
at most two of the maximal intervals in $f(n)\setminus g(n)$ --- one, $I_n$,
at the top and and one, $J_n$, at the bottom.
The two sequences $\langle i_n\rangle_{n\in\omega}$ 
(bottom points of the~$I_n$)
and $\langle j_n\rangle_{n\in\omega}$ 
(top points of the~$J_n$) 
determine cut points $i_u$ and $j_u$ of~$\II_u$, which cannot belong to~$L$.
Therefore we can enlarge $m$ and $\alpha$ such that 
$\bigl\{n:i_n,j_n\notin[a_{m,n},b_{\alpha,n}]\bigr\}$ is in~$u$.
A convexity argument will now establish that 
$\bigl\{n:\bigl(f(n)\setminus g(n)\bigr)\cap[a_{m,n},b_{\alpha,n}]
  =\emptyset\bigr\}$ belongs to~$u$.
The same argument, interchanging $f$ and $g$ will yield our final $m$ 
and $\alpha$.
\end{proof}

\begin{lemma}\label{A-h-cover-L}
$L\subset\cl{A_f}$ if and only if there are $m$ and $\alpha<\omega_1$ 
such that the set
$\bigl\{n:[a_{m,n},b_{\alpha,n}]\subseteq f(n)\bigr\}$ belongs to~$u$.
\end{lemma}
\begin{proof}
Apply Lemma~\ref{lemma:Af=Ag} to $f$ and the constant function $n\mapsto\II$.
\end{proof}

Now we are ready to perform the construction of~$g_\omega$, given
subsets $\{f_k\}_{k\le\omega}$ and $\{g_k\}_{k<\omega}$ 
of $\powerinfront{\omega}{\cR}$ such that the map
$\cl{A_{f_k}}\cap L \mapsto \cl{A_{g_k}}\cap[L,1_u]$ ($k<\omega$)
satisfies the conditions (1)--(4) from our list.

The conditions that need to be met are
\begin{itemize} 
\item[(a)] $L\cap\cl A_{f_\omega}=L\cap\cl A_{g_\omega}$; 
\item[(b)] if $L\subseteq \cl A_{f_k}\cup \cl A_{f_\omega}$
      then $[L,1_u]\subseteq \cl A_{g_k}\cup \cl A_{g_\omega}$; and
\item[(c)] if $F\subseteq\omega$ is finite and
      $L\cap\cl A_{f_\omega}\cap\bigcap_{l\in F}\cl A_{f_l}=\emptyset$ then
      $[L,1_u]\cap\cl A_{g_\omega}\cap\bigcap_{l\in F}\cl A_{g_l}=\emptyset$. 
\end{itemize}
The first condition takes care of (1) and (4) in our list, except possibly
when $\cl A_{f_\omega}\cap L=\emptyset$ but in that case it suffices
to let $g_\omega$ be the constant function $n\mapsto\emptyset$.
The second and third condition ensure (2) and (3) respectively.
There is one more condition that we need to keep the recursion alive;
it is needed to take care of combinations of (b) and (c):
if $L\subseteq \cl A_{f_k}\cup \cl A_{f_\omega}$ and 
$L\cap\cl A_{f_\omega}\cap\bigcap_{l\in F}\cl A_{f_l}=\emptyset$
then we must have room to be able to ensure that both
$[L,1_u]\subseteq \cl A_{g_k}\cup \cl A_{g_\omega}$ and 
$[L,1_u]\cap\cl A_{g_\omega}\cap\bigcap_{l\in F}\cl A_{g_l}=\emptyset$.
Note that the antecedent implies that, in the subspace~$L$, 
the intersection $L\cap\bigcap_{l\in F}\cl A_{f_l}$ is contained in 
the interior of $L\cap\cl A_{f_k}$.
A moment's reflection shows that we need
\begin{itemize}
\item[(d)] if $L\cap\bigcap_{l\in F}\cl A_{f_l}$ is contained in 
           $\interior_L L\cap\cl A_{f_k}$ then
           $[L,1_u]\cap\bigcap_{l\in F}\cl A_{g_l}$ is contained in 
           $\interior_{[L,1_u]} [L,1_u]\cap\cl A_{f_k}$.
\end{itemize}
For every $k$ as in (b) choose $m_k$ and $\alpha_k$ as per 
Lemma~\ref{A-h-cover-L} such that 
$U_k=
\bigl\{n:[a_{m_k,n},b_{\alpha_k,n}]\subseteq f_k(n)\cup f_\omega(n)\bigr\}$
belongs to~$u$.
Likewise, for every $F$ as in~(c) choose $m_F$ and $\alpha_F$ as per
Lemma~\ref{lemma:AfcapL=0} such that
$U_F=\bigl\{n:[a_{m_F,n},b_{\alpha_F,n}]\cap f_\omega(n)\cap 
      \bigcap_{l\in F} f_l(n)=\emptyset\bigr\}$
belongs to~$u$.
And, finally, for every pair $(F,k)$ as in (d) 
(with $F$ finite but with $k\le\omega$ in this case) 
choose $m_{F,k}$ 
and $\alpha_{F,k}$, and $U_{F,k}\in u$ such that for every $n\in U_{F,k}$
we have 
$[a_{m_{F,k},n},b_{\alpha_{F,k},n}]\cap\bigcap_{l\in F}f_l(n)
   \subseteq\interior f_k(n)$ 
and
$[a_{m_{F,k},n},1]\cap\bigcap_{l\in F}g_l(n)\subseteq\interior g_k(n)$
(the latter only if $k<\omega$ of course).

We fix an ordinal~$\alpha$ larger than the $\alpha_k$, 
$\alpha_F$ and $\alpha_{F,k}$ by~$\alpha$ and use it instead in the 
definitions of the sets~$U_k$, $U_F$ and~$U_{F,k}$ --- 
they will still belong to~$u$.
Next take a decreasing sequence $\langle V_p\rangle_{p\in\omega}$
of elements of~$u$ such that $V_p$ is a subset of
\begin{itemize}
\item $U_k$ whenever $k<p$;
\item $U_F$ whenever $F\subseteq p$; and 
\item $U_{F,k}$ whenever $F\subseteq p$ and $k<p$ or $k=\omega$. 
\end{itemize}
In addition we can, and will, assume that whenever $F\subseteq p$
and $L\cap\bigcap_{l\in F}\cl A_f=\emptyset$ then
$[b_\alpha,1]\cap\bigcap_{l\in F}g_l(n)=\emptyset$ --- that this is possible
follows from the assumption that (c) holds for $\max F$.

Now were are truly ready to define~$g_\omega$.
If $n\notin V_0$ define $g_\omega(n)=\II$.
In case $n\in V_p\setminus V_{p+1}$ observe first that if $k<p$ is as in~(b)
and $F\subseteq p$ is as in~(c) then $(F,k)$ is as in~(d) 
so that certainly 
$$
[a_{m_{F,k}},1]\cap\bigcap_{l\in F}g_l(n)\subseteq\interior g_k(n).\eqno(*)
$$
Define $g_\omega(n)$ as the union of $f_\omega(n)\cap[0,b_\alpha(n)]$
and an element $h(n)$ of~$\cR$ that is a subset of~$[b_\alpha(n),1]$
and satisfies
\begin{itemize}
\item $h(n)\cup g_k(n)\supseteq[b_\alpha(n),1]$ whenever $k<p$ is as in (b);
\item $h(n)\cap\bigcap_{l\in F}g_l(n)=\emptyset$ whenever $F\subseteq p$
      is as in~(c); and
\item $h(n)\supseteq [b_{\alpha,n},1]\cap\bigcap_{l\in F}g_l(n)$
       whenever $(F,\omega)$ is as in~(d).
\end{itemize}
This is possible because of $(*)$ and because 
$\bigcap_{l\in F}g_l(n) \cap\bigcap_{l\in G}g_l(n)=\emptyset$
whenever $F$ is as in~(c) and $(G,\omega)$~is as in~(d).
This gives us just enough room to choose~$h(n)$.

It is now routine to verify that all conditions on~$g_\omega$ are
met $u$-often:
e.g., if $F\subseteq\omega$ is finite and
      $L\cap\cl A_{f_\omega}\cap\bigcap_{l\in F}\cl A_{f_l}=\emptyset$
then $[a_{m_F,n},1]\cap g_\omega(n)\cap\bigcap_{l\in F}g_l(n)=\emptyset$
for all $n\in V_p$, where $p=1+\max F$.

\subsection{Further considerations}

The proof in the previous section can be used to show that, under~$\CH$,
all other layers of the continuum~$\II_u$ are retracts of~$\II_u$.
If the layer is a point then this is clear.
If the layer~$L$ is non-trivial then the cofinality of $[0_u,L)$ and
the coinitiality of $(L,1_u]$ are $\omega_1$.
It is then a matter of making the proof of Theorem~\ref{retraction} symmetric
to get our retraction $r:\II_u\to L$.
The details can be found in~\cite{vanderSteeg2003}.

The fixed-point free homeomorphism $h:L\to L$ from 
Theorem~\ref{thm:all.indec.not.fpp} can then be used to construct 
another witness to $s^*(\II_u)\neq0$, 
almost exactly as in the proof of Theorem~\ref{thm:sssI-u.neq.0}.

\section{Remarks}

The results of this paper grew out of an attempt to find non-metric
counterexamples to Lelek's conjecture.
The fairly easy proof, indicated after 
Corollary~\ref{cor:sigmastar.Hstar.nonzero},
that $\Hstar$ is not chainable, which also works for layers of countable
cofinality lead us to consider $\II_u$ as a possible candidate.

A secondary goal was to convert any non-metric counterexample into a
metric one by an application of the L\"owenheim-Skolem theorem
(\cite{Hodges1997}*{Section~3.1})
to its lattice of closed sets.
This produces a countable sublattice with exactly the same
(first-order) lattice-theoretic properties; its Wallman representation
space, see~\cite{W}, is a metrizable continuum with many properties
in common with the starting space, e.g., covering dimension
unicoherence, (hereditary) indecomposability, \dots, 
see~\cite{vanderSteeg2003}*{Chapter 2} for a comprehensive list. 

The results of this paper cast doubt of the possibility of adding
(non-)chainability and span (non)zero (of any kind) to this list.
The reason for this is that the family 
$\cR_u=\{\cl A_f\cap\II_u:f\in\powerinfront\omega\cR\}$ is isomorphic
to the ultrapower of~$\cR$ (from the proof Theorem~\ref{retraction})
by the ultrafilter~$u$;
this follows in essence from the equivalence of 
$\cl A_f\cap\II_u=\cl A_g\cap\II_u$ and
$\{n:f(n)=g(n)\}\in u$.
By the \L os Ultraproduct Theorem
(\cite{Hodges1997}*{Theorem 8.5.3})
we see that $\cR$ and $\cR_u$ have the same first-order lattice theoretic
properties yet their Wallman representations, $\II$ and $\II_u$ respectively,
differ in chainability and in various kinds of span 
(all kinds if $\CH$ is assumed).

Chainability is a property that can be read off from a lattice base for the 
closed sets (or dually for the open sets): using compactness one readily shows
that a continuum is chainable iff every basic open cover has a chain refinement
from the base.
Thus we deduce that chainability is not a first-order property of the
lattice base.

For span (non)zero there are two possibilities:
it cannot be read off from a base or, if it can be, it is not a first-order
property of the lattice base.

\section{Questions}

The remarks in the previous section suggest lots of questions.
We mention the more important ones.

\begin{question}
Is there a non-metric counterexample to any one version of Lelek's
conjecture?  
\end{question}

It should be noted that, as mentioned in~\cite{Da}, H. Cook has shown
that the dyadic solenoid has symmetric span zero.

In spite of the results on $\II$ and $\II_u$ it is still possible
that the L\"owenheim-Skolem method may convert a non-metric counterexample
into a metric one.
The reason for this is that $\cR_u$ is special base for the closed sets
of~$\II_u$ and not an elementary sublattice of its lattice of closed sets.

\begin{question}
If $L$ is an elementary sublattice of the full lattice of closed sets of
the continuum~$X$, does its Wallman representation
inherit (non-)chainability and or span (non)zero from~$X$?  
\end{question}

Section 3.7 of \cite{vanderSteeg2003} gives a positive answer for very
special sublattices, but unfortunately except for span zero.
Further, more specialized, questions can be found in that reference.

The corollaries in Section~\ref{sec:surj-spans} were derived from 
Theorem~\ref{retraction}, which needed $\CH$ in its proof.
This clearly suggests the question whether a more insightful analysis
of the structure of the $\II_u$ and the use of more intricate combinatorics
will make the use of~$\CH$ unnecessary.

\begin{question}
Can one show in $\ZFC$ only that all spans of $\Hstar$ and $\II_u$
are nonzero? 
\end{question}

It would already be of interest if one could find at least one~$u$
such that all spans of~$\II_u$ are nonzero.

We have shown implicitly that the fixed-point property like chainability
and span zero in that $\II$ has it but $\II_u$ does not, at least
under $\CH$.

\begin{question}
Is there in $\ZFC$ at least one $u$ such that $\II_u$ does not have the
fixed-point property?  
\end{question}



\end{document}